# A FORMULA FOR FINDING A POTENTIAL FROM NODAL LINES

JOYCE R. MCLAUGHLIN AND OLE H. HALD

ABSTRACT. In this announcement we consider an eigenvalue problem which arises in the study of rectangular membranes. The mathematical model is an elliptic equation, in potential form, with Dirichlet boundary conditions. We have shown that the potential is uniquely determined, up to an additive constant, by a subset of the nodal lines of the eigenfunctions. A formula is given which, when the additive constant is fixed, yields an approximation to the potential at a dense set of points. An estimate is presented for the error made by the formula.

## INTRODUCTION

In this summary we consider the Dirichlet eigenvalue problem for the operator $(-\triangle + q)$ on a rectangle. This problem arises in the study of rectangular, vibrating membranes. The goal here is to solve the inverse problem: find $q$ from the nodal line positions of the eigenfunctions. We show that, for almost all rectangles and for sufficiently smooth $q$, a potential $q$ whose integral average is zero can be uniquely determined from a subset of the nodal lines of the eigenfunctions.

The theorems, which we present in this paper, extend the one-dimensional results of McLaughlin and Hald [8], [4], [5], where inverse nodal problems are defined. There the authors showed that the one-dimensional potential for the Sturm-Liouville problem with Dirichlet or mixed boundary conditions is uniquely determined, up to an additive constant, by a dense set of nodes. The analysis of the Sturm-Liouville problem depends on perturbation theory. The proofs of the perturbation results are reasonably straightforward. This is possible because the difference between consecutive eigenvalues increases as the order of the eigenvalues increases. The perturbation results provide asymptotic forms for large eigenvalues, the corresponding eigenfunctions and nodal positions. The uniqueness theorem follows.

For the two-dimensional problem we also require asymptotic forms for the eigenvalues and eigenfunctions and approximate location of nodal lines. The difference between the one-dimensional and the higher dimensional cases is that in the higher dimensional case the eigenvalues are not well separated. Even for those rectangles for which all eigenvalues are distinct, arbitrarily large eigenvalues can be arbitrarily

Received by the editors October 6, 1993, and, in revised form, September 9, 1994.

1991 *Mathematics Subject Classification.* Primary 35R30, 35P20, 73D50.

Partial support for the first author's research came from ONR grant N00014-91J-1166 and from NSF grant VPW-8902967. Partial support for the second author's research came from NSF grant DMS-9003033.

The research announced here was presented at the 1994 International Congress of Mathematicians in Zürich, Switzerland.

close together. These small differences can become small divisors in an asymptotic expansion; see Kato [7]. To handle this, we start with the $q = 0$ case and show there is a dense, well-defined set of rectangles for which the eigenvalues are well separated. The rectangles are selected so that the square of the ratio of the sides is not well approximated by rational numbers. Our condition is almost the same as that given by Moser [9]. For each acceptable rectangle, we choose "good" eigenvalue, eigenfunction pairs that have the following properties. These eigenvalues are at least a given distance from their nearest neighbors and a large distance from a selected set of neighbors. The corresponding eigenfunctions have nearly the same number of oscillations in the $x$ and $y$ directions. Using number theoretic and geometric results, it can be shown that almost all eigenvalue, eigenfunction pairs satisfy these properties.

Using the $q = 0$ case as a base problem, the application of perturbation theory yields leading order terms for "good" eigenvalue, eigenfunction pairs for the $q \neq 0$ case. Explicit estimates are obtained for remainder terms. Our technique is strongly influenced by perturbation results for the two- and three-dimensional periodic problem by Feldman, Knörrer, and Trubowitz [3] and Friedlander [2]. Instead of $L^2$ estimates, we need and derive $L^\infty$ estimates for the difference between the eigenfunctions and their leading order terms; we require the $L^\infty$ bound so that we can determine the approximate location of the nodal lines for the $q \neq 0$ case.

Having established asymptotic forms for almost all the eigenfunctions and having established the $L^\infty$ bound for the difference between the eigenfunctions and their leading order terms, we encounter the second major difficulty. A nodal domain for an eigenfunction for the $q \neq 0$ case is not necessarily a small perturbation of a nodal domain for the corresponding eigenfunction for the $q = 0$ case. This follows from the fact that the nodal line intersections for the $q = 0$ case are not usually present in the nodal line patterns for the $q \neq 0$ case. The fact that nodal line intersections are not stable under perturbation has been shown by Uhlenbeck; see [11], [12]. Figure 1 on page 245 shows typical nodal domains for the corresponding eigenfunctions for the $q = 0$ case and the $q \neq 0$ case.

In our analysis we select approximate nodal domains which are subsets of nodal domains for the $q \neq 0$ problem. The approximate nodal domains are small perturbations of nodal domains for the $q = 0$ problem. The eigenfunction for the $q \neq 0$ case is either zero or nearly zero on the boundary of the approximate nodal domain.

In all of the above, we address properties of the forward problem. Only then can we address the inverse problem.

We solve the inverse nodal problem. A formula is given for approximating $q$ at a dense set of points using only a subset of the nodal line positions. We show that, at each point, $q$ can be approximated by the difference of two eigenvalues. One eigenvalue is an eigenvalue for the $q = 0$ Dirichlet problem defined on the rectangle. The other eigenvalue is the first eigenvalue for the $q = 0$ Dirichlet problem defined on an approximate nodal domain. An estimate is given for the error in making the approximation. A uniqueness theorem follows. The proofs, which are presented in [6], use the Rayleigh-Ritz variational formulation of the eigenvalue problem and the perturbation results described above.

In the body of this paper we present statements of our key results. Here we do not give the sharpest possible results but give a presentation that makes the ideas more accessible to the reader. The complete set of results together with the proofs are contained in [6].

## Main results

We begin with the problem

$$
\begin{aligned}
-\Delta u + qu &= \lambda u\,, & (x\,,\,y)\,&\in\,R\,, \\
u &= 0\,, & (x\,,\,y)\,&\in\,\partial R\,,
\end{aligned}
\tag{1}
$$

where $R = [0\,,\,\pi/a] \times [0\,,\,\pi]\,$, $a > 1$ and $q \in C_0^\infty(R)\,$. Our goal is to solve the inverse problem: Find $q$ from the positions of the nodal lines. To attain this goal, we must first study the forward problem. When $q = 0\,$, that is, when we have the problem,

$$
\begin{aligned}
-\Delta u &= \lambda u\,, & (x\,,\,y)\,&\in\,R\,, \\
u &= 0\,, & (x\,,\,y)\,&\in\,\partial R\,,
\end{aligned}
$$

we can label the eigenvalue, eigenfunction pairs as

$$
\lambda_{\alpha 0} = |\alpha|^2 = |(an\,,\,m)|^2 = a^2 n^2 + m^2\,,
$$

$$
u_{\alpha 0} = \frac{2\sqrt{a}}{\pi} \sin anx \sin my\,,
$$

respectively. Here $\alpha$ is in the two-dimensional lattice

$$
L(a) = \{\alpha = (an\,,\,m)|\,n\,,\,m = 1\,,\,2\,,\,3\,,\,\ldots\}\,.
$$

We must first establish a perturbation result for the perturbation from $q = 0$ to $q \neq 0\,$. We can do this for almost all of the eigenvalues $\{\lambda_{\alpha 0}\}_{\alpha \in L}$ when the parameter $a\,$, which defines the rectangle $R\,$, is well chosen. Specifically we choose $a^2$ so that it is not well approximated by rational numbers; i.e. $a$ is in the set

$$
V = \Bigg\{ a > 1\,|\ \text{there exists}\ \ 0 < \delta < \frac{1}{24}\ \ \text{and}\ \ k > 0\,,
$$

$$
\text{such that for all p,}\ \ q > 0 : |a^2 - \frac{p}{q}| > \frac{k}{q^{2+\delta}} \Bigg\}\,.
$$

It can be shown that in any subinterval, $J = (1\,,\,a_0)\,$, meas $(J \backslash V) = 0\,$. Further we can write down specific irrational numbers which are in $V\,$. For example we can show that $\sqrt{e-1} \in V\,$; and it follows from a theorem of Roth [10] that if $a^2$ is an irrational, algebraic number, then $a \in V\,$.

Our perturbation result is

**Theorem 1.** *Let $a$ belong to $V\,$. Let $q \in C_0^\infty(R)$ have mean value $0\,$. There exists an exceptional set of lattice points $M(a) \subset L(a)$ depending on $a$ and an upper bound on the derivatives of $q$ of order no more than $100$ such that*

a) $M(a)$ *has density* 0 *in* $L(a)$ *in the sense that*

$$\lim_{r \to \infty} \frac{\#\{\alpha \in L(a)\backslash M(a) \mid |\alpha| < r\}}{\#\{\alpha \in L(a) \mid |\alpha| < r\}} = 1.$$

b) *For every* $\alpha \in L(a)\backslash M(a)$, *there is a unique eigenvalue of the variable coefficient problem* (1) *which we denote as* $\lambda_{\alpha q}$ *satisfying*

$$|\lambda_{\alpha q} - \lambda_{\alpha 0}| \le |\alpha|^{-15/8}.$$

*A suitable multiple of the corresponding eigenfunction satisfies*

$$\left\| u_{\alpha q} - u_{\alpha 0} - \sum_{\beta \ne \alpha} \frac{(q u_{\alpha 0}, u_{\beta 0})}{|\alpha|^2 - |\beta|^2} u_{\beta 0} \right\|_\infty \le \sqrt{a} |\alpha|^{-15/8}.$$

*Remark.* Our techniques for the proof of Theorem 1 are strongly influenced by [1], [2], [3], [7], [13].

*Remark.* We note that there are two properties of the $L^\infty$ bound for the eigenfunctions that are important to us. One is that we choose "good" $\alpha$ so that we have

$$\left\| \sum_{\beta \ne \alpha} \frac{(q u_{\alpha 0}, u_{\beta 0})}{|\alpha|^2 - |\beta|^2} u_{\beta 0} \right\|_\infty < \sqrt{a} |\alpha|^{-15/16}.$$

We use this to establish the position of the nodal lines for $u_{\alpha q}$. The second property is that near the intersections of the nodal lines of $u_{\alpha 0}$ the bound

$$|u_{\alpha q}| < \sqrt{a} |\alpha|^{-15/8}$$

holds; this is used to establish the error bounds for the formula for the potential $q$.

We have established the $L^\infty$ bound for $u_{\alpha q} - u_{\alpha 0}$ (see the Remark above) so that we can determine the approximate location of nodal lines and hence define approximate nodal domains for $u_{\alpha q}$. We first exhibit a nodal domain for $u_{\alpha 0}$, $\alpha = (an, m)$; it is

$$\Omega_0 = \left\{ (x, y) \left| \left| x - \frac{(n_1 + \frac{1}{2})\pi}{an} \right| < \frac{\pi}{2an}, \left| y - \frac{(m_1 + \frac{1}{2})\pi}{m} \right| < \frac{\pi}{2m} \right. \right\},$$

for some $n_1 = 0, 1, ..., n - 1$, $m_1 = 0, 1, ..., m - 1$, and note that a nodal domain of $u_{\alpha q}$ with $\alpha \in L\backslash M$ is a connected component of $\{(x,y)|u_{\alpha q} \ne 0\}$. Now we have a difficulty because a nodal domain of $u_{\alpha q}$ may be a large perturbation of a nodal domain for $u_{\alpha 0}$. Figure 1 shows the nodal domains for $u_{\alpha 0}$ and $u_{\alpha q}$ for a typical case.

Note that the nodal line intersections for $u_{\alpha 0}$ are not present in the nodal line pattern for $u_{\alpha q}$. See [11], [12] for the proof that nodal line intersections are unstable under perturbation.

In Figure 1, the "natural" perturbation of $\Omega_0$ is the much larger nodal domain $\Omega$. We will select a subdomain of $\Omega$ as our approximate nodal domain. To do this more precisely, let $d = 1/[(3.9)|\alpha|^{15/16}]$, $d_1 = d/(an)$, $d_2 = d/m$, and let nodal domains for $u_{\alpha 0}$                    nodal domains for $u_{\alpha q}$

FIGURE 1

$\Omega_1 \subset \Omega_0 \subset \Omega_2$ be

$$\Omega_1 = \left\{ (x\,,\,y) \,\middle|\, \left| x - \frac{(n_1 + \frac{1}{2})\pi}{an} \right| < \frac{\pi}{2an} - d_1\,, \quad \left| y - \frac{(m_1 + \frac{1}{2})\pi}{m} \right| < \frac{\pi}{2m} - d_2 \right\}\,,$$

$$\Omega_2 = \left\{ (x\,,\,y) \,\middle|\, \left| x - \frac{(n_1 + \frac{1}{2})\pi}{an} \right| < \frac{\pi}{2an} + d_1\,, \quad \left| y - \frac{(m_1 + \frac{1}{2})\pi}{m} \right| < \frac{\pi}{2m} + d_2 \right\}\,.$$

Figure 2 exhibits the three domains $\Omega_0$, $\Omega_1$, $\Omega_2$.

We can then show that $u_{\alpha q} \neq 0$ in $\Omega_1$. This allows us to choose a natural nodal domain $\Omega$ associated with $\Omega_0$. It is the nodal domain of $u_{\alpha q}$ which contains $\Omega_1$. Our approximate nodal domain is a subset of $\Omega$ and is as follows.

**Definition.** Let $\alpha \in L(a) \backslash M(a)$ and let $\Omega_0$ be a nodal domain for $u_{\alpha 0}$. Define the corresponding $\Omega_1 \subset \Omega_0 \subset \Omega_2$ as above. Let $\Omega$ be the nodal domain of $u_{\alpha q}$ which contains $\Omega_1$. Then the perturbation $\Omega'$ of $\Omega_0$ is the connected component of $\Omega \cap \Omega_2$ which contains $\Omega_1$.

Figure 3 depicts a choice for $\Omega'$.

These are the main theorems of this paper.

**Theorem 2.** *Let* $q \in C_0^\infty(R)$ *with mean value zero, and let* $\alpha \in L \backslash M$. *Let* $\Omega_0$ *be a nodal domain of* $u_{\alpha 0}$, *and let* $\Omega'$ *be the approximate nodal domain of* $u_{\alpha q}$

FIGURE 2

FIGURE 3. Construction of $\Omega'$.

*defined above. Let* $\lambda_{1,0}(\Omega')$ *be the first eigenvalue of*

$$-\Delta u = \lambda u, \qquad (x, y) \in \Omega',$$
$$u = 0, \qquad (x, y) \in \partial\Omega'.$$

*Then, there exists* $(\bar{x}, \bar{y}) \in \Omega'$ *with*

$$|q(\bar{x}, \bar{y}) - [|\alpha|^2 - \lambda_{1,0}(\Omega')]| < |\alpha|^{-7/4}.$$

Because $M$ has density 0 in $L$, the set of points $(\bar{x}, \bar{y})$ for all nodal domains of $u_{\alpha 0}, \alpha \in L \backslash M$, is dense in $R$. This proves

**Theorem 3.** *Let* $q \in C_0^\infty(R)$ *with mean value zero and let* $\alpha \in L \backslash M$. *Then* $q$ *is uniquely determined by the nodal lines of the eigenfunctions for*

$$(-\Delta + q)u = \lambda u, \qquad (x, y) \in R,$$
$$u = 0, \qquad (x, y) \in \partial R.$$

*Remark.* We can improve on Theorem 3 showing that $q$ is uniquely determined by a subset of the nodal lines. To do this we use the property that, when $\alpha \in L \backslash M$, $an$ and $m$ are comparable. From this it follows that $q$ is uniquely determined by the nodal lines of $u_{\alpha^k q}$ for any set $\alpha^k \in L \backslash M$, $k = 1, 2, ...$ where $\mid \alpha^k \mid \to \infty$ as $k \to \infty$.

## Acknowledgments

The authors have benefited from discussions with E. Trubowitz, L. C. Evans, H. Knörrer, K. Ribet, A. van der Poorten, P. Vojta and G. Bergman. The authors wish to thank the referee for her/his suggestions for the revision of this announcement.

Rensselaer Polytechnic Institute, Troy, New York 12181
*E-mail address*: mclauj@rpi.edu

University of California, Berkeley, California 94720